\documentclass[a4paper,11pt]{amsart}

\usepackage[dvips]{graphicx}

\usepackage{ulem, fancybox}

\begin{document}

\title{Deciding nonconstructibility of $3$-balls \\
with strict spanning edges}
\author{Satoshi KAMEI}
\date{}

\maketitle

\begin{center}{School of Computer Science, 

Tokyo University of Technology, 

1404-1 Katakura, Hachioji-shi Tokyo 192-0982, Japan

e-mail: kamei@cs.teu.ac.jp}\end{center}

\begin{abstract}
In this paper, we consider constructibility of simplicial $3$-balls. In many cases, examining $1$-dimensional subcomplexes of a simplicial $3$-ball is efficient to solve the decision problem whether the simplicial $3$-ball is constructible or not. From the point of view, we consider the case where a simplicial $3$-ball has spanning edges and present a sufficient condition for nonconstructibility.
\end{abstract}

\vspace*{\baselineskip}

\begin{center}{\small Key words: simplicial complex; shellability; constructibility}\end{center}

\vspace*{\baselineskip}

\section{Introduction}

Shellability and constructibility are combinatorial concepts relating to decompositions of simplicial balls. Since each shellable simplicial ball is constructible, the two concepts form a hierarchy. Therefore to check constructibility of a simplicial ball is useful to decide whether the simplicial ball is shellable or not.

In~\cite{H1}, Hachimori presented a necessary and sufficient condition for constructibility of a reduced $3$-ball which has no interior vertices. Further he proved that the condition is still necessary and sufficient in the case where a reduced $3$-ball has at most two interior vertices in~\cite{H2}. In~\cite{K}, the author researched the case where a $3$-ball has interior vertices and no edge which connects interior vertices, and presented a sufficient condition for nonconstructibility. In this paper, we deal with a more general case.

We prepare some notions. See~\cite{B1} for the terminology of simplicial complexes. Further see~\cite{Z1} and~\cite{Z2} for the definition of shellability. 
First, we recall the definition of constructibility of pure simplicial complexes.

\vspace{\baselineskip}
\noindent
{\bf Definition.} A pure $d$-dimensional simplicial complex $C$ is {\it constructible} if

\noindent
(1) $C$ is a simplex, or

\noindent
(2) there exist $d$-dimensional constructible subcomplexes $C_1$ and 
$C_2$ such that $C = C_1 \cup C_2$ and that $C_1 \cap C_2$ is a 
$(d-1)$-dimensional constructible complex.

\vspace{\baselineskip}

A simplicial $3$-ball $B$ is {\it reduced} if each $2$-face of $B \setminus \partial B$ has at most one edge which is contained in $\partial B$. As shown in ~\cite{H1}, there are operations for reduction and we can obtain a reduced $3$-ball from any simplicial $3$-ball in finite steps of the operations. 
Furthermore, if some reduced $3$-ball $B_2$ is obtained from a $3$-ball $B_1$ by the reduction and $B_2$ is nonconstructible, then the $3$-ball $B_1$ is also nonconstructible. Thus it is essential to consider a sufficient condition for nonconstructibility of reduced $3$-balls.


An edge of a simplicial $3$-ball $B$ is called {\it spanning} if the edge is not contained in $\partial B$ and its end vertices are contained in $\partial B$. 
In~\cite{K}, we define the notion ``strict" only for a spanning edge of a simplicial $3$-ball which has interior vertices and no edge the end vertices of which are contained in the interior of the $3$-ball. We redefine the notion to deal with more general cases. For the purpose, we make preparations.

\vspace*{\baselineskip}

\noindent
{\bf Definition.} Let $B$ be a simplicial $3$-ball. A vertex contained in the interior of $B$ is called an {\it interior vertex} of $B$. An edge is called an {\it interior edge} of $B$ if the both end of the edge are contained in the interior of $B$. A closed $1$-dimensional simplicial complex which consists of interior vertices and interior edges of $B$ is called an {\it interior graph} of $B$. A connected component of an interior graph is called an {\it interior graph component}. If an interior graph component is maximal in $B$, the interior graph component is called a {\it maximal interior graph component} of $B$. 

\vspace*{\baselineskip}

Notice that an edge is not called an interior edge if the edge has an end vertex which is contained in $\partial B$. 

In the followings, we denote by ${\rm Star}_B I$ the closed star neighbourhood of $I$ with respect to $B$. 

\vspace*{\baselineskip}

\noindent
{\bf Definition.} A spanning edge $e$ is called {\it strict} if there is no maximal interior graph component $I$ of $B$ such that both end vertices of $e$ are contained in exactly one connected component of $\partial B \cap {\rm Star}_B I$. 

\vspace*{\baselineskip}

Now we state the main theorem.

\vspace*{\baselineskip}

\noindent
{\bf Theorem 2.3.} {\it If a reduced $3$-ball has spanning edges and all of the spanning edges are strict, then the $3$-ball is nonconstructible.}

\vspace*{\baselineskip}

In Section 2, we prove Theorem 2.3. In Section 3, we present two examples.

\section{Main argument}

First, we prepare technical lemmas for the proof of Theorem 2.3.

\vspace*{\baselineskip}

\noindent
{\bf Lemma 2.1.} Let $D$ be a simplicial $2$-ball and $V$ be the set which consists of all vertices of $D$. Let $W$ be a set of some vertices which are contained in $\partial D$. We denote by $I$ the closed $1$-dimensional maximal subcomplex of $D$ the vertices of which coincide with $V \setminus W$. We assume that $C$ is connected and that any edge which connects vertices of $W$ is contained in $\partial D$. Then ${\rm Star}_D I$ coincides with $D$.

\vspace*{\baselineskip}

\noindent
{\it Proof.} There is no $2$-simplex in $D$ all vertices of which are contained in $W$, if otherwise there exists an edge which is not contained in $\partial D$ and connects vertices of $W$. Thus, for each vertex in $W$, there exists an edge which connects the vertex and a vertex in $V \setminus W$. Therefore all $2$-simplices in $D$ are contained in ${\rm Star}_D I$, thus ${\rm Star}_D I$ coincides with $D$. \qed

\vspace*{\baselineskip}

\noindent
{\bf Lemma 2.2.} We follow the symbols $D$, $V$, $W$ and $I$ in Lemma 2.1. Assume that $I$ has more than one components. Then there exists an edge which is not contained in $\partial D$ and connects vertices of $W$.

\begin{figure}[htbp]
 \begin{center}
  \includegraphics{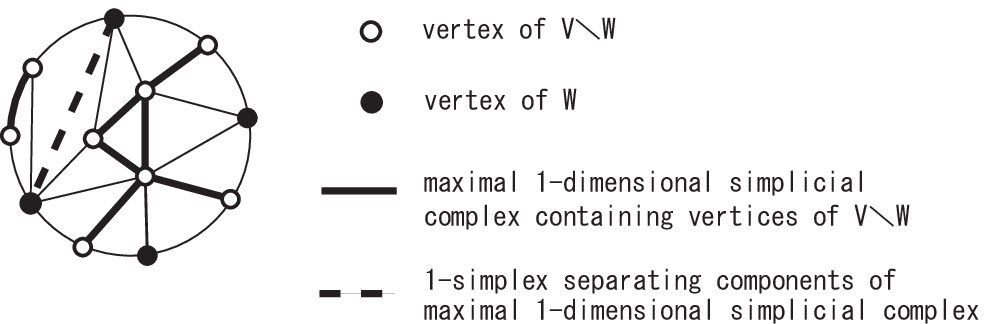}
 \end{center}
   \caption{}

\end{figure}

\vspace*{\baselineskip}

\noindent
{\it Proof.} Let $I_1$ be a component of $I$. Then ${\rm Star}_D I_1$ dose not coincide with $D$ and all vertices of ${\rm Star}_D I_1 \cap \overline{ \{D \setminus {\rm Star}_D I_1 \} }$ are contained in $W$. Thus there exists an edge which satisfies the conditions. \qed

\vspace*{\baselineskip}

Further, we prepare some notions for the proof of Theorem 2.3.

\vspace*{\baselineskip}

\noindent
{\bf Definition.} Let $D$ be a simplicial $2$-ball which is a subcomplex of a simplicial $3$-ball. Assume that $D$ contains spanning edges of the $3$-ball. A spanning 
edge $e$ is {\it outermost} if $e$ cuts off a simplicial $2$-ball $\Delta$ from $D$ which contains no spanning edge except $e$. The simplicial $2$-ball $\Delta$ is called an {\it outermost disk} of $D$.



\vspace*{\baselineskip}

\noindent
{\bf Definition.} Let $B$ be a constructible $3$-ball and $B'$ be a subcomplex of $B$ such that $B' =B'_1 \cup B'_2$ is part of the construction of $B$. Then the $2$-ball $B'_1 \cap B'_2$ is called a {\it divide}.

\vspace*{\baselineskip}

Notice that the interior of a divide is contained in the interior of $B$. Therefore the interior of a divide contains no vertices which are contained in $\partial B$. 

Now we prove the main theorem.

\vspace*{\baselineskip}

\noindent
{\bf Theorem 2.3.} {\it If a reduced $3$-ball has spanning edges and all of the spanning edges are strict, then the $3$-ball is nonconstructible.}

\vspace*{\baselineskip}

\noindent
{\it Proof.} Let $B$ be a reduced $3$-ball satisfying the hypothesis. We assume that $B$ is constructible and consider a deconstruction of $B$ as the followings. First, we set $B^{(1)} = B$. At each step of the deconstruction, a $3$-ball $B^{(i)}$ is decomposed into two $3$-balls $B^{(i)}_1$ and $B^{(i)}_2$ such that $B^{(i)} = B^{(i)}_1 \cup B^{(i)}_2$ is part of the construction of $B$. Since $B$ is constructible, there exists such a sequence of decompositions and $B$ is decomposed into simplices in finite steps. Let $D^{(i)} = B^{(i)}_1 \cap B^{(i)}_2$, thus $D^{(i)}$ is the divide which appears in the boundaries of decomposed $3$-balls at the $i$-th step of the deconstruction.

We choose the step $i_0$ so that for the first time throughout the deconstruction there appear spanning edges on the boundaries of decomposed $3$-balls. Thus the spanning edges are contained in $D^{(i_0)}$. 
We choose an outermost spanning edge $e$ which cuts off an outermost disk $\Delta$ from $D^{(i_0)}$. 
Let $\cup I_n$ be the union of all maximal interior graph components of $B$. Then $\Delta \cap \{ \cup I_n \}$ is connected, if otherwise there would be a spanning edge except $e$ in $\Delta$ from Lemma 2.2 and it contradicts the assumption that $\Delta$ is outermost. Let $\Gamma^{(i_0)}$ be $\Delta \cap \{ \cup I_n \}$ and $P^{(i_0)}$ be $\overline{\partial \Delta \setminus e}$. From Lemma 2.1, ${\rm Star}_{\Delta} \Gamma^{(i_0)}$ coincides with $\Delta$, thus the path $P^{(i_0)}$ is contained in ${\rm Star}_{B^{(i_0)}} \Gamma^{(i_0)}$.

If all vertices of $P^{(i_0)}$ are contained in $\partial B$, the spanning edge $e$ is not strict and it contradicts the assumption that all spanning edges are strict. Thus there exist vertices of $P^{(i_0)}$ which are not contained in $\partial B$.  Let $v_1$ and $v_2$ be the end vertices of $e$. Notice that $v_1$ and $v_2$ are contained in $\partial B$ because $e$ is a spanning edge. In the followings, we construct paths $P^{(j)}$ which connect $v_1$ and $v_2$ and which are contained in $\partial B^{(j)}$ for $1 \le j \le i_0-1$ inductively. 


\begin{figure}[htbp]
 \begin{center}
  \includegraphics{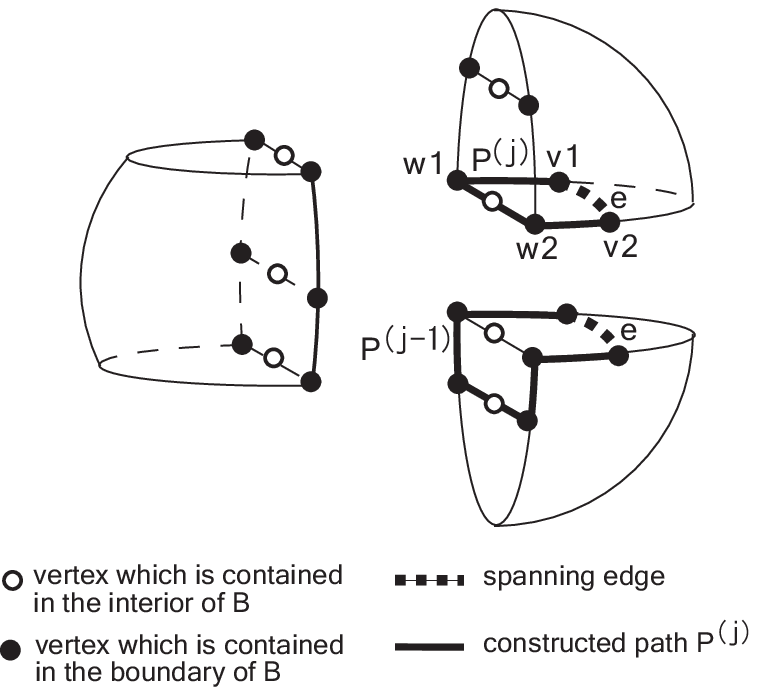}
 \end{center}
   \caption{}

\end{figure}

\vspace*{\baselineskip}


We start with $P^{(i_0)}$ and $\Gamma^{(i_0)}$. If $P^{(i_0)} \cap \{D^{(i_0-1)} \setminus \partial D^{(i_0-1)}\}$ is empty, we set $P^{(i_0-1)} = P^{(i_0)}$ and $\Gamma^{(i_0-1)} = \Gamma^{(i_0)}$. Otherwise, let $w_1$ and $w_2$ be the vertices of $P^{(i_0)} \cap \partial D^{(i_0-1)}$ which are the nearest to $v_1$ and $v_2$ respectively. If $P^{(i_0)} \cap \partial D^{(i_0-1)}$ contains $v_1$ and $v_2$, we set $w_j = v_j$ for $j = 1, 2$. Notice that there exists more than one vertices which are contained in $P^{(i_0)} \cap \partial D^{(i_0-1)}$ because $P^{(i_0)} \cap \{D^{(i_0-1)} \setminus \partial D^{(i_0-1)}\}$ is not empty and $\partial D^{(i_0-1)}$ does not separate $v_1$ and $v_2$ on $\partial B^{(i_0-1)}$. We switch the subpath of $P^{(i_0)}$ with a connected subcomplex of $\partial D^{(i_0-1)}$ both of which connect $w_1$ and $w_2$. Although there are two connected subcomplexes of $\partial D^{(i_0-1)}$ the end vertices of which are $w_1$ and $w_2$, we may choose either of them. We denote by $P^{(i_0-1)}$ the constructed path the end vertices of which are $v_1$ and $v_2$. 
Further, we construct an interior graph component $\Gamma^{(i_0-1)}$. Since there is no spanning edge which is contained in $D^{(i_0-1)}$, $D^{(i_0-1)} \cap \{ \cup I_n \}$ is connected. It is obvious that there exists no $2$-simplex of $D^{(i_0-1)}$ all vertices of which are contained in $\partial B$. Thus, $\partial D^{(i_0-1)}$ is contained in ${\rm Star}_{D^{(i_0-1)}} D^{(i_0-1)} \cap \{ \cup I_n \}$ from Lemma 2.1. Since there exists a vertex of $P^{(i_0)}$ which is not contained in $\partial B$, $\{ D^{(i_0-1)} \cap \{\cup I_n\} \} \cup \Gamma^{(i_0)}$ is connected. We denote it by $\Gamma^{(i_0-1)}$. Then, the path $P^{(i_0-1)}$ is contained in ${\rm Star}_{B^{(i_0-1)}} \Gamma^{(i_0-1)}$.

Continuously, we construct $P^{(j)}$ and $\Gamma^{(j)}$ from $j =i_0-1$ to $1$. Then, we obtain the path $P^{(1)}$ which is contained in $\partial B$ and ${\rm Star}_B \Gamma^{(1)}$. Thus the spanning edge $e$ is not strict and it contradicts our assumption. \qed

\section{Examples}

In this section, we construct two examples. The first example is a $3$-ball
which satisfies the hypothesis of Theorem 2.3.

\begin{figure}[htbp]
 \begin{center}
  \includegraphics{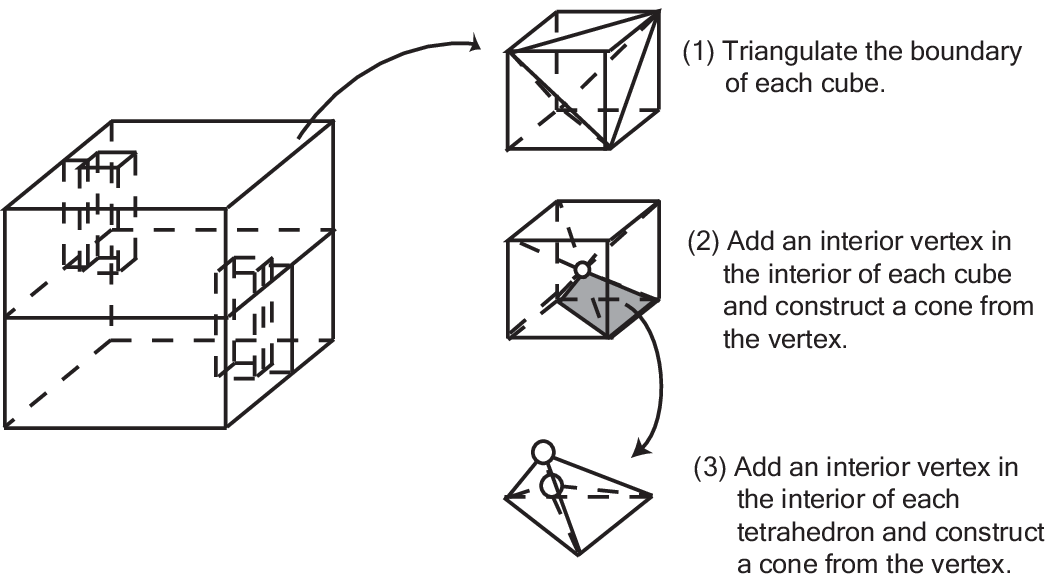}
 \end{center}
   \caption{}

\end{figure}

\vspace*{\baselineskip}

\noindent
{\bf Example 3.1.} We construct an example based on Example 4.1 in~\cite{K}. 
First, we recall the example. Consider the $3$-ball which is depicted in the left side of Figure 3. The walls of the ball are made of one layer of cubes. There are two tunnels one of which connects the upper space and 
the lower floor, and the other of which connects the lower space and the upper 
floor. 

We triangulate the $3$-ball as the following. First, triangulate the $2$-skeleton of the $3$-ball. Next, add an interior vertex in each cube, and construct a cone from the vertex over the triangulated boundary of each cube. Notice that there is an appropriate triangulation of the $2$-skeleton of the $3$-ball so that the triangulated $3$-ball is reduced. See~\cite{H1} and~\cite{K} for more details. The constructed simplicial $3$-ball satisfies the hypothesis of Theorem 3.1 in~\cite{K}. Thus the simplicial $3$-ball is not constructible.

Now we add an interior vertex in each tetrahedron and construct a cone from the vertex over the boundary of each tetrahedron. There exists exactly one maximal interior graph component in each cube. The constructed simplicial $3$-ball satisfies the hypothesis of Theorem 2.3, thus the $3$-ball is still nonconstructible.

\vspace*{\baselineskip}

\vspace*{\baselineskip}
The second example is a shellable $3$-ball which has a strict 
spanning edge and nonstrict spanning edges. 

\vspace*{\baselineskip}

\begin{figure}[htbp]
 \begin{center}
  \includegraphics{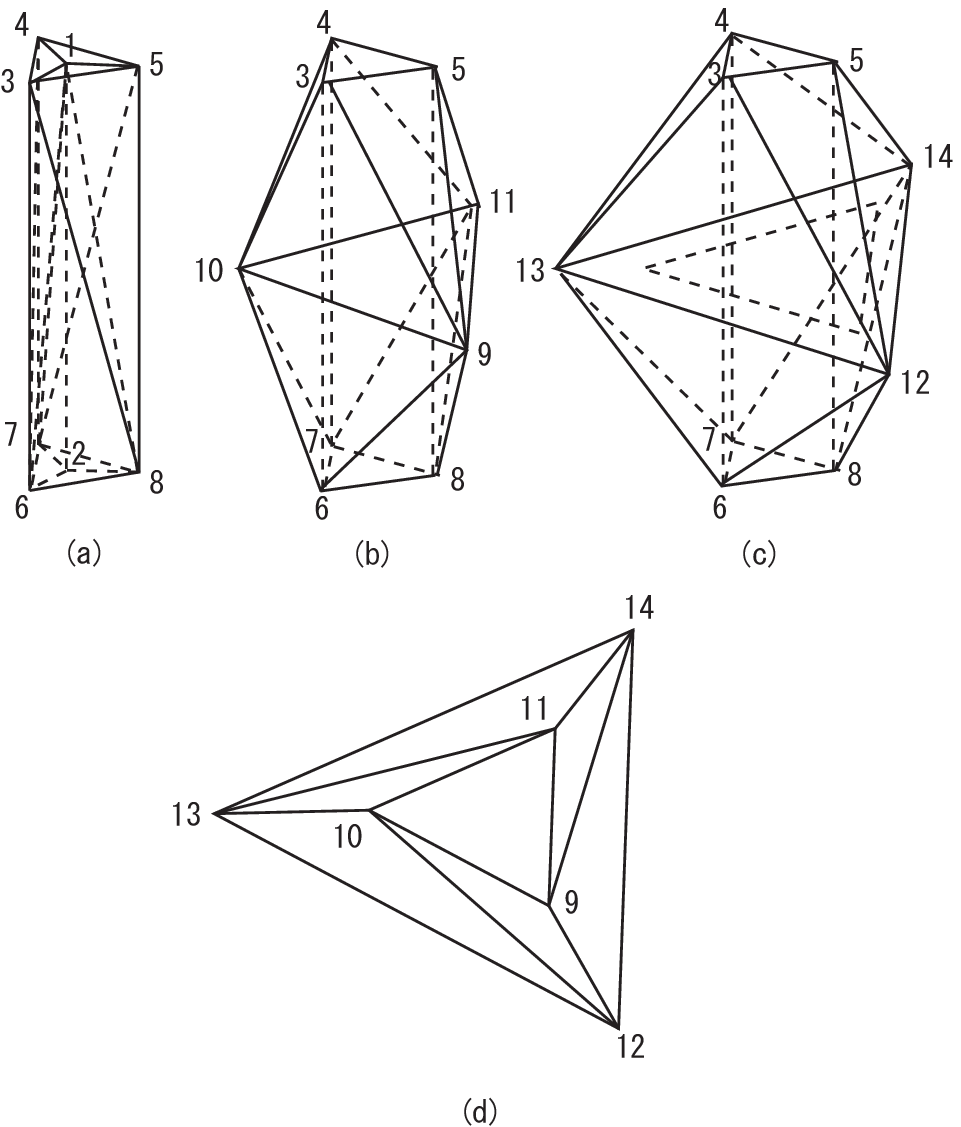}
 \end{center}
   \caption{}

\end{figure}

\noindent
{\bf Example 3.2.} This example is constructed as indicated in Figure 4. First, we consider a solid trigonal pillar ``3-4-5-6-7-8". Add edges ``3-8",``5-7" and ``4-6", then the boundary of the pillar is triangulated. Add vertices ``1" and ``2" in the interior of $2$-simplices ``3-4-5" and ``6-7-8" respectively. Connect the vertex ``1" and the vertices ``4, 5, 6", and the vertex``2" and ``7, 8, 9" by edges respectively. Add a spanning edge ``1-2" and triangulate the pillar as Figure 4(a). Put a simplicial $1$-sphere ``9-10-11" such that the pillar goes through the $1$-sphere. Connect the point ``9" and the points ``3, 5, 6, 8" by edges respectively. Also connect the point ``10" and the points ``3, 4, 6, 7", the point ``11"  and the points `` 4, 5, 7, 8" in the same way. Fill tetrahedra which have the edges added at the previous step and the edges on the boundary of the pillar. Further put a simplicial $1$-sphere ``12-13-14" as Figure 4(d). Connect the point ``12" and the points ``3, 5, 6, 8, 9, 10", the point ``13" and the points ``3, 4, 6, 7, 10, 11", and the point ``14" and the points ``4, 5, 7, 8, 9, 11" by edges respectively as Figure 4(c) and Figure 4(d). Fill tetrahedra the same as the previous step. 

The constructed $3$-ball is reduced and has exactly one maximal interior graph component "9-10-11". Further the $3$-ball has one strict spanning edge ``1-2" and three nonstrict spanning edges ``3-6", ``4-7", ``5-8". We can easily check that the $3$-ball is shellable thus it is also constructible. This example suggests that we cannot relax the latter hypothesis of Theorem 2.3, that is, all spanning edges of a reduced $3$-ball are strict.

\vspace*{\baselineskip}

The following is the list of all facets of this example.

\[ \begin{array}{lllll}
	  \{1,2,6,7 \} &  \{1,2,7,8\}  & \{1,2,6,8\} & \{1,3,5,8\}  & \{1,3,6,8\} \\
	   \{1,3,4,6\}  & \{1,4,6,7\} & \{1,4,5,7\}  & \{1,5,7,8\} \\
	 	  \{5,8,9,11\} & \{3,5,8,9\} & \{3,6,8,9\} & \{3,6,9,10\} & \{3,4,6,10\} \\
	 	  \{4,6,7,10\} & \{4,7,10,11\} & \{4,5,7,11\} & \{5,7,8,11\} & \\
	 	  \{3,4,10,13\} & \{3,10,12,13\} & \{3,9,10,12\} & \{3,5,9,12\} & \{5,9,12,14\} \\	 	  
	  \{5,9,11,14\} &  \{4,5,11,14\} &  \{4,11,13,14\} & \{4,10,11,13\} & \\
	  \{6,7,10,13\} &  \{6,10,12,13\} &  \{6,9,10,12\} & \{6,8,9,12\} & \{8,9,12,14\} \\
 \{8,9,11,14\} & \{7,8,11,14\} & \{7,11,13,14\} &  \{7,10,11,13\} \\
\end{array} \]

\vspace*{\baselineskip}

\noindent
{\bf Acknowledgement.} The author would like to thank Dr. Masahiro Hachimori for his helpful advice and discussions for the construction of Example 3.2.

\pagebreak

\end{document}